\documentclass[12pt]{article}
\usepackage{amsmath,amssymb,amsfonts}
\usepackage[margin=1in]{geometry}
\usepackage[english,russian,shorthands=off]{babel}
\usepackage{tikz,tikz-cd,graphicx}
\def\Point#1#2#3{\node [fill=black,inner sep=1pt,label=#1:{#2}] at (#3) {}}

\newcount\pno \def\mn{\medskip\noindent}
\def\makenote#1#2{\expandafter\def\csname#1\endcsname{\advance\pno by 1
\medskip\noindent\textbf{#2 \the\pno. }}}\makenote{defn}{Definition}
\makenote{ex}{Example}\makenote{thm}{Theorem}\makenote{lem}{Lemma}
\makenote{rem}{Remark}\makenote{prop}{Proposition}\makenote{cor}{Corollary}
\def\avnt#1#2#3{{\setbox0=\hbox{$#1{#2#3}{\int}$ }
\vcenter{\hbox{$#2#3$ }}\kern-.6\wd0}}\def\mod#1{\ (\mathrm{mod}\ #1)}
\def\avint{\mathchoice{\avnt\displaystyle\textstyle-}%
{\avnt\textstyle\scriptstyle-}{\avnt\scriptstyle\scriptscriptstyle-}%
{\avnt\scriptscriptstyle\scriptscriptstyle-}\!\int}
\def\pf{\medskip\noindent\textit{Proof.} }\def\qed{\hfill$\square$}
\def\F{\mathbb F}\def\P{\mathcal P}\def\R{\mathbb R}\def\Z{\mathbb Z}
\def\RP{\mathbb{RP}}\def\d{\mathrm d}
\def\Sym{\mathrm{Sym}}\def\Gr{\widetilde{\mathrm{Gr}}}
\def\<{\langle}\def\>{\rangle}\def\D{\nabla}\def\eps{\varepsilon}
\def\iff{\Leftrightarrow}
\begin{document}
\begin{center}{\Large\textbf{Note on Lebesgue's universal cover problem}}
\\Yizhen Chen

\bigskip\textbf{Abstract}

\medskip\begin{minipage}{.75\textwidth}
A universal cover is a convex set that covers all sets of diameter 1
after some rotation, reflection, and translation.
First, we show that a regular dodecahedron that circumscribes a ball
of diameter 1 is not a universal cover in $\R^3$, answering a question
of Chakerian [K]. Second, we improve the bounds
on the minimum volume of a universal cover in $\R^3$ to
$(0.545193,0.655984)$. Third, we give an
infinite family of centrally symmetric polytopes in $\R^d$ with
$d(d+1)$ facets that circumscribe the unit ball and are universal covers.
This is the maximum possible number and proves a conjecture of Makeev [M2].
Finally, on the least number of facets of a centrally symmetric polytope
in $\R^d$ that circumscribe the unit ball and is not a universal cover,
we generalize Makeev's lower bound [M1] and improve his upper bound [M3].
\end{minipage}

\bigskip\textsc{1. Introduction}
\end{center}

\noindent\textbf{Definition 1. }\pno=1
A convex set $K\subset\R^d$ is a \textbf{universal cover} if for every
set $X\subset\R^d$ of diameter 1, there exists $A\in O(d)$ and
$b\in\R^d$ such that $X\subset A(K)+b$.

Let $K^d(v_1,\ldots,v_n):=
\{x\in\R^d:\forall i,|\<x,v_i\>|\le\frac12\|v_i\|\}$. When this notation
is used, we assume that $v_i\in\R^d$ are nonzero vectors
such that $v_i\ne\pm v_j$ for $i\ne j$ and
$\mathrm{span}(v_1,\ldots,v_n)=\R^d$, making this set a centrally
symmetric polytope with $2n$ facets.

\medskip In 1914, H.~Lebesgue proposed the problem of finding
universal covers with the least volume. For example, the unit cube is
a universal cover in all dimensions $d$, and it is optimal when $d=1$.
When $d=2$, the best bounds for the minimum volume of a universal cover
are $(0.83222,0.8440936)$ obtained by Gibbs [Gi],
and Brass and Sharifi [BS].

By P\'al's theorem [MMO, Thm.~7.2.2], every set of diameter 1 is contained
in a set of constant width 1, so it suffices to consider only such $X$.

When $d=3$, Makeev gave the current best upper bound
$\frac12(7-4\sqrt2)\approx0.671573$
by showing that the intersection of the rhombic dodecahedron
$R=K^3(e_i\pm e_j:i\ne j)$ and $\{x\in\R^3:\<x,e_i\>\le\frac12\}$
is a universal cover [M5]. He conjectured that $K^d(\{e_i-e_j\in V:i\ne j\})$,
where $V=\{x\in\R^{d+1}:\sum_{i=0}^dx_i=0\}$,
is a universal cover for any dimension $d>1$.
Kuperberg gave another proof that $R$ is a universal cover,
and conjectured that the regular cross polytope (in $\R^4$) and
the regular dodecahedron (in $\R^3$) that circumscribe a ball of diameter 1
are universal covers [K]. The last of these questions was asked by
Don Chakerian.

In this note, we shall answer Chakerian's question in Section 2, and give new
bounds for the minimum volume of a universal cover in $\R^3$
in Section 3.

\thm (i) The regular dodecahedron in $\R^3$ that circumscribes a
sphere of diameter 1 is not a universal cover.
\\(ii) The dual of the $3_{21}$ polytope in $\R^7$ that circumscribes
a sphere of diameter 1 is not a universal cover. The polytope in $\R^{23}$
that circumscribes a sphere of diameter 1, whose facets are normal to
an orbit of 552 vectors of the Conway group $\mathrm{Co}_3$,
is also not a universal cover.

\thm In $\R^3$, the minimum volume of a universal cover is in the interval
$(0.545193,0.655984)$.

\medskip For arbitrary $d>0$, Makeev showed that $K^d(v_1,\ldots,v_n)$
is not a universal cover if $n>\frac12d(d+1)$, and he conjectured that
this bound is sharp [M2]. He found an infinite
family of universal covers of the form $K^d(v_1,\ldots,v_n)$
where $n=\frac14d(d+4)$ for even $d$ and $n=\frac14(d^2+4d-1)$ for odd
$d>1$ [M4]. In Section 4, we shall prove the conjecture of Makeev,
giving an infinite family of universal covers
of the form $K^d(v_1,\ldots,v_n)$ where $n=\frac12d(d+1)$.

\thm Let $P_{2j}\subset S^1\subset\R^2$ be a polygon with
$2j$ vertices, one of them being $(1,0)$, and with central symmetry
and the reflection symmetry $(x,y)\mapsto(x,-y)$.
Let $\pm v_1,\ldots,\pm v_n$ ($n=\frac12d(d+1)$) be the vertices of
$\prod_{j=1}^kP_{2(4j-1)}\subset\R^{2k}$ for $d=2k$, and
$[-1,1]\times\prod_{j=1}^kP_{2(4j+1)}\subset\R^{2k+1}$ for $d=2k+1$.
Then $K^d(v_1,\ldots,v_n)$ is a universal cover.

\medskip Makeev also studied the question of when there exists
at least one polytope $K=K^d(v_1,\ldots,v_n)$ that is not a universal cover. 
He showed that any such $K$ with $n=d+1$ is a universal cover, and
gave better lower bounds (for the least possible $n$)
when $d\in\{4,8\}$ [M1]. On the other hand,
he found a family of such $K$ with $n=3d-2$ that is not a universal cover [M3].
We shall generalize the lower bound, which we used for proving Theorem 3,
to arbitrary dimensions. We shall also improve the upper bound
to $n=2d$. In particular, this is sharp when $d$ is a power of 2.

\prop Let $K=K^d(v_1,\ldots,v_n)$, where $d>2$.
\\(i) Any $K$ with $n=d+2$ and a rotational symmetry
taking $v_{d+1}$ to $-v_{d+1}$ and $v_{d+2}$ to $-v_{d+2}$
is a universal cover. In particular, if $d$ is even, this is
always satisfied by $-1$.
\\(ii) Any $K$ with $n=d+k+1$ is a universal cover, where $d=2^em$,
$m$ is odd, and \[ k=
\begin{cases}2&e=1,\\2^e-1&m>1\ne e,\\2^e-2&m=1.\end{cases}\]

\thm Let $D=K^3(e_1,e_2,e_3,e_1+e_2,e_2+e_3,e_3+e_1)$, and
$P=K^2(e_1,e_2,e_1\pm e_2)$ (regular octagon), then $P^k$ is not a universal
cover when $d=2k$, and $D\times P^{k-1}$ is not a universal cover
when $d=2k+1$.

\begin{center}\textsc{2. Specific Non-Universal Covers}\end{center}

We use the following notations for convenience.

\defn The \textbf{support function} of a set $X\subset\R^d$ is
$h_X:\R^d\to\R$, $v\mapsto\sup_{x\in X}\<x,v\>$.
The \textbf{width} and shift of
$X$ in direction $v\in S^{d-1}$ are $w_X(v)=h_X(v)+h_X(-v)$ and
$s_X(v)=\frac12(h_X(v)-h_X(-v))$ respectively.

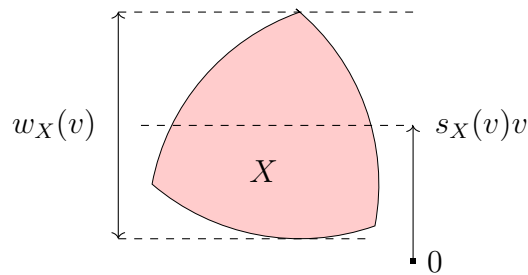
\begin{center}\begin{tikzpicture}[scale=3]
\path [draw=black, fill=red!20, line width=.1mm] (0,1.5)
arc [radius=1, start angle=109.38, end angle=169.38]
arc [radius=1, start angle=-130.62, end angle=-70.62]
arc [radius=1, start angle=-10.62, end angle=50.62] -- cycle;
\draw [dashed] (-0.8,0.5) -- (0.3,0.5); \draw [dashed] (-0.7,1) -- (0.5,1);
\draw [dashed] (-0.8,1.5) -- (0.5,1.5); \draw [<->] (-0.8,0.5) -- (-0.8,1.5);
\draw [->] (0.5,0.4) -- (0.5,1);
\node [label=right:{$s_X(v)v$}] at (0.5,1) {};
\node [label=left:{$w_X(v)$}] at (-0.8,1) {};
\node [label=left:{$X$}] at (0,0.8) {};
\node [fill=black,inner sep=1pt,label=right:{0}] at (0.5,0.4) {};
\end{tikzpicture}\end{center}

\noindent\advance\pno by 1\textbf{Lemma \the\pno.}
Given any odd smooth function $s_0:S^{d-1}\to\R$, there exists
a set $X\subset\R^d$ of constant width 1 with $s_X=\eps s_0$ for
sufficiently small $\eps>0$.

\pf Let $h(v)=\|v\|(1+\eps s_0(\frac v{\|v\|}))$,
then $\D^2h(v)=\|v\|^{-1}(I-\frac{v\otimes v}{\|v\|^2})+
\eps\D^2h_0(v)$ where $h_0(v)=\|v\|s_0(\frac v{\|v\|})$.
Note that $h(\lambda v)=\lambda h(v)$ implies
$\D h(\lambda v)=\D h(v)$, so $\D^2h(v)v=\D_v\D h(v)=0$.
On the sphere, $I-v\otimes v$ is the projection to $TS^{d-1}$,
so its eigenvalues are 1. Thus for $\eps>0$ sufficiently small,
$\D^2h>0$ on $TS^{n-1}$, so $h$ is convex. Hence there exists a
set $X\subset\R^d$ of constant width 1 with support function $h$ by
[Gr, Thm.~4.3].\qed

\mn\textit{Proof of Theorem 2(i).} Let $H_3$ be the space of degree 3
\textbf{spherical harmonics}---homogeneous polynomials $f$ in 3 variables
such that $\Delta f=0$. Let $s_0\in H_3$. By Lemma 8, there is a set
$X\subset\R^3$ of constant width 1 such that $s_X=\eps s_0$
for some $\eps>0$.

Let $K=K^3(v_1,\ldots,v_6)$ be a regular dodecahedron.
Suppose $K$ is a universal cover, then in particular
there exist $A\in O(3)$ and $b\in\R^3$ such that $s_X(Av_i)=\<b,Av_i\>$
for $1\le i\le6$. Then $(A^\top s_0)(v_i)=\<\eps^{-1}A^\top b,v_i\>$.
The evaluation map ev from odd functions on $S^2$ to odd functions on
$\{\pm v_i:1\le i\le6\}$ is $A_5$-equivariant, where $A_5$ acts
as the rotational symmetries of $K$.
The two sides of the equation are in ev$(H_3)$ and ev$(V)$,
where $V$ is the space of linear functions $\R^3\to\R$.

\textit{Claim:} $H_3$ and $V$ do not share
irreducible components as $A_5$-representations,
so $s_0(Av_i)=0$ for $1\le i\le6$.

We first prove the theorem assuming this claim. The equation
$(A^\top s_0)(v_i)=\<\eps^{-1}A^\top b,v_i\>$ now implies $s_0(Av_i)=0$
for $1\le i\le6$. In cylindrical
coordinates $x=(r\cos\theta,r\sin\theta,z)$, let $s_0(x)=3r^2z-5z^3$,
which equals $3z-5z^3$ on the unit sphere. Since $s_0(Av_i)=0$,
for fixed unit vector $v=A^\top e_3$, $\<Av_i,e_3\>^2\in\{0,\frac35\}$.
On the other hand, we have $\sum_{i=1}^6v_i\<v_i,\cdot\>=2\cdot\mathrm{id}$,
so $\sum_{i=1}^6\<Av_i,e_3\>^2=\sum_{i=1}^6\<v_i,A^\top e_3\>^2
=2\|A^\top e_3\|^2=2$. This cannot be the sum of several $\frac35$.

\textit{Proof of the Claim.} The rotational symmetry group $A_5$
of the regular dodecahedron consists of: 1 identity, 15 rotations by $\pi$
around lines connecting midpoints of opposite edges, $2\times10$
rotations by $\frac{2\pi}3$ around lines connecting opposite vertices,
and $2\times6$ rotations by $\frac{2k\pi}5$ around lines connecting
centers of opposite faces ($k=1,2$). These are all the conjugacy classes
of $A_5$, and we shall denote them by their rotation angles.

We compute the character of the representation $V$ of $A_5$
as follows. In cylindrical coordinates,
we have $V=\mathrm{span}\{re^{i\theta},re^{-i\theta},z\}$
(e.g.\ the function $z$ is $\<e_3,\cdot\>$).
Every rotation in $SO(3)$ is conjugate to a rotation
$\theta\mapsto\theta+\alpha$ around the $z$-axis,
which, in our basis of $V$, has trace
$\chi_V(\alpha)=e^{i\alpha}+e^{-i\alpha}+1=1+2\cos\alpha$. Restricting
to $A_5\subset SO(3)$, we find that $V$ is isomorphic to the
irreducible representation $V_3$ of $A_5$.

\begin{center}\begin{tabular}{c|ccccc}
class&1&$\pi$&$\frac{2\pi}3$&$\frac{2\pi}5$&$\frac{4\pi}5$\\\hline
size&1&15&20&12&12\\\hline
$\chi_1$&1&1&1&1&1\\
$\chi_4$&4&0&1&$-1$&$-1$\\
$\chi_5$&5&1&$-1$&0&0\\
$\chi_3$&3&$-1$&0&$\phi_+$&$\phi_-$\\
$\chi_3'$&3&$-1$&0&$\phi_-$&$\phi_+$\\\hline
$\chi_V$&3&$-1$&0&$\phi_+$&$\phi_-$\\
$\chi_{H_3}$&7&$-1$&1&$\phi_--1$&$\phi_+-1$
\end{tabular}\end{center}

On the space $\Sym^3(V)$ of homogeneous degree 3 polynomials
in $x,y,z$ (where $A_5\subset SO(3)$ acts), the image of the Laplacian
$\Delta$ is $V$ and the kernel is $H_3\cong\mathrm{Sym}^3(V)\ominus V
\cong V_3'\oplus V_4$. So $H_3$ and $V$ do not share irreducible
components. In the table above, $\chi_i$ is the character of $V_i$,
and $\phi_\pm=\frac12(1\pm\sqrt5)$.\qed

\medskip The key property making this proof work is the Claim:
$G$-representations $\Sym^3(V)$ and $V=\R^d$ ($d=3$ in this case)
do not share irreducible components, where $G$ is a finite subgroup
of $O(d)$ and the normal vectors $\pm v_1,\ldots,\pm v_n$ of
facets of the polytope (the regular dodecahedron in this case)
are a $G$-orbit. Equivalently, Goethals and Seidel showed that
the $G$-orbit is a \textbf{spherical 5-design}, i.e. \begin{equation}
\frac1{2n}\sum_{i=1}^n(p(v_i)+p(-v_i))=\avint_{S^{d-1}}p \end{equation}
for all polynomials $p:\R^d\to\R$ of degree $\le5$
[GS, Thm.~6.1]. Delsarte, Goethals, and
Seidel showed that a spherical 5-design in $\R^d$
has at least $d(d+1)$ vectors [DGS, Thm.~5.12], while Makeev showed
that any universal cover of the form $K^d(v_1,\ldots,v_n)$ has
$n\le\frac12d(d+1)$ [M2, Thm.~2], so we look for centrally symmetric
spherical 5-designs with exactly $d(d+1)$ vectors. For $d>3$, only
two are known: the 56 vertices of the $3_{21}$ polytope in $\R^7$,
and the 552 vectors in $\R^{23}$ that form an orbit of the Conway group
$\mathrm{Co}_3$ [DGS, Ex.~8.3]. It turns out that the polytopes
$K^d(v_1,\ldots,v_n)$ formed from these vectors are not universal
covers either, although the argument in the proof of Theorem 2(i) fails.
Numerical computation shows that random choices of degree 3 spherical
harmonic as $s_0$ often work, which motivates our argument, a refinement
of Makeev's method in [M2] (see Proposition 13 below).

\mn\textit{Proof of Theorem 2(ii).} Let $\pm v_i$
($1\le i\le n=\frac12d(d+1)$) be the normal vectors of facets of
the polytope. Let $\pi:\R^d\to\R^{d-2}$
be the projection to the last $d-2$ coordinates. We shall show that
there is an odd polynomial $f:\R^{d-2}\to\R$ such that no $A\in SO(d)$
satisfy $(\pi^*f)(Av_i)=\<b,v_i\>$ for all $i$. Then the polytope does
not cover the set $X$ of constant width 1 with $s_X=\eps f$ given
by Lemma 8.

Let $L=\{(\<b,v_i\>)\in\R^n:b\in\R^d\}$. Let $W_A\subset\R^n$ be
the vector space of ``potential values'' of $((\pi^*f)(Av_i))$ for
odd functions $f$ (all $\R^2\subset\R^d$ below are embedded as the
first two coordinates): \begin{align*}W_A:=\{(w_i)\in\R^n:
w_i=0&\text{ if }\pi(Av_i)=0\text{ (iff }v_i\in A^{-1}(\R^2)),\\
w_i=\pm w_j&\text{ if }\pi(Av_i)=\pm\pi(Av_j)\text{ (iff }
v_i\pm v_j\in A^{-1}(\R^2))\}.\end{align*}
Let $\P$ be the set of odd polynomials $\R^{d-2}\to\R$ of degree $<d(d+1)$,
then the map $\P\to W_A$, $f\mapsto((\pi^*f)(Av_i))$, is surjective:
let \begin{equation}\{x_1,\ldots,x_{2m}\}=\{\pi(Av_i):1\le i\le n\}
\text{ where }x_{m+i}=-x_m,\end{equation}
then $((\pi^*f)(Av_i))$ is determined by the values $f(x_i)$. Let
$f_i(x)=\prod_{j\ne i}\frac{\<x-x_j,x_i-x_j\>}{\|x_i-x_j\|^2}$,
then $f_i(x)-f_i(-x)$ is an odd polynomial of degree $<2m\le n$
that takes $x_j$ to $\delta_{ij}$.

We have a smooth fiber bundle $X=SO(d)/SO(2)\to\Gr_2(\R^d)$,
$A\mapsto Z:=A^{-1}(\R^2)$, with fibers $SO(d-2)$. There are only finitely
many possible $W_A$, which depends only on $Z$:
if \[k:=\dim\mathrm{span}(Z\cap\{v_i,v_i\pm v_j:i\ne j\}) \] is 0,
then $W_A=\R^n$; $k=1$ happens only in the finitely many
$(d-2)$-dimensional smooth submanifolds (such as $\{Z:v_i\in Z\}$)
of $\Gr_2(\R^d)$; $k=2$ happens only at finitely many points
$Z\in\Gr_2(\R^d)$. Outside these closed submanifolds, $k=0$.

Over each one of these submanifolds $Y$ (including the open one with $k=0$),
$W_A$ is constant and the map $\phi:\P\times X|_Y\to(W_A+L)/L$,
$(f,A)\mapsto((\pi^*f)(Av_i))$, is a submersion by the previous paragraph.
So $\phi^{-1}(0)$ is a closed submanifold of dimension
\begin{align*}\dim\phi^{-1}(0)&=\dim\P+\dim X|_Y-\dim((W_A+L)/L)\\
&=\dim\P+\dim Y+\dim SO(d-2)-\dim W_A+\dim(W_A\cap L).\end{align*}
It suffices to show that $\dim\phi^{-1}(0)<\dim\P$. Then $\phi^{-1}(0)$
projects to a measure zero set in $\P$ corresponding to those polynomials
$f$ without the desired property. For the last term $W_A\cap L$,
$(\<b,v_i\>)\in L$ is in $W_A$ iff $\<b,v_i\>=0$ whenever $v_i\in Z$
and $\<b,v_i\pm v_j\>=0$ whenever $v_i\pm v_j\in Z$,
iff $b\in(\mathrm{span}(Z\cap\{v_i,v_i\pm v_j:i\ne j\}))^\perp$, so
$\dim(W_A\cap L)=d-k$. Hence it suffices to show that
\[\dim W_A>\frac12d(d-1)-(k-1)(d-1)\] for $k\in\{1,2\}$ (it holds when $k=0$).

\mn\textit{The bound for $k=2$.}

Now we use the property that $\{\pm v_1,\ldots,\pm v_n\}$ is a
spherical 5-design. The map $S^1\times B^{d-2}\to S^{d-1}$,
$(z,w)\mapsto(z\sqrt{1-\|w\|^2},w)$ (where $B^{d-2}\subset\R^{d-2}$
is the closed unit ball), has Jacobian determinant 1,
so for an even polynomial $f:\R^{d-2}\to\R$ of degree 4,
\begin{align*} \avint_{S^{d-1}}\pi^*f&=\avint_{B^{d-2}}\avint_{S^1}
(\pi^*f)(z\sqrt{1-\|w\|^2},w)\,\d z\,\d w=\avint_{B^{d-2}}f\\
\text{and }\avint_{S^{d-1}}\pi^*f&=\frac1n\sum_{i=1}^n(\pi^*f)(Av_i)
=\sum_{i=0}^ma_if(x_i)\text{ where }x_0:=0\end{align*}
by (1) and (2), where $a_i>0$ for $i>0$, and $a_0\ge0$.
Noskov and Schmid showed that this implies $m>\frac12(d-1)(d-2)$
for $d\ge6$ not an odd square, which is the case here ($d\in\{7,23\}$) [NS].

\mn\textit{The bound for $k=1$.}

The vertices of the $3_{21}$ polytope are
$\pm v_{ij}/\sqrt{24}\in u^\perp\subset\R^8$
where $v_{ij}=u-4(e_i+e_j)$, $u=\sum_{i=1}^8e_i$, and $1\le i<j\le8$ [C].
The following (and their coordinate permutations) are all the possible cases
for $W_A$. So $\dim W_A\ge n-6=\frac12d(d-1)+1$.
\begin{center}\begin{tabular}{llll}\hline
$\in A^{-1}(\R^2)$&in $W_A$&$\in A^{-1}(\R^2)$&in $W_A$\\\hline
$v_{12}$&$w_{12}=0$&$v_{13}-v_{23}$&$w_{1i}=w_{2i}$ ($3\le i\le8$)\\
$v_{12}-v_{34}$&$w_{12}=w_{34}$&$v_{13}+v_{23}$&$w_{13}+w_{23}=0$\\
$v_{12}+v_{34}$&\multicolumn{3}{l}{%
$w_{12}+w_{34}=w_{13}+w_{24}=w_{14}+w_{23}=$}\\&\multicolumn{3}{l}{%
$w_{56}+w_{78}=w_{57}+w_{68}=w_{58}+w_{67}=0$}\\\hline
\end{tabular}\end{center}

For a description of the vectors $v_1,\ldots,v_n$ in the $d=23$ case,
as a projection of a subset $S$ of the Leech lattice, see
[CS, Section 10.3.6]. The set $S=\{u_1,\ldots,u_n\}\subset\R^{24}$ contains
certain coordinate permutations of $\sqrt2(e_1+e_2)$ and
$(\sum_{i=1}^8e_i)/\sqrt2$. Let $x=(4e_1+\sum_{i=1}^{24}e_i)/\sqrt8$, then
$\|x\|^2=6$, $\|u_i\|^2=4$, and the vectors $v_i$ are in the hyperplane
$\<u_i,x\>=3$. Also $\<u_i,u_j\>\in\{1,2\}$ for $i\ne j$ as
the Leech lattice has no vectors of norm $\sqrt2$. The vectors
$v_i=u_i-\frac x2$ are the projections of $u_i$ to $x^\perp\cong\R^{23}$.

From the properties above, we have $\|v_i\|^2=\frac52$ and
$\<v_i,v_j\>=\pm\frac12$ for $i\ne j$. Suppose
$v_1\pm v_2+a(v_3\pm v_4)=0$ (each $\pm$ may be different),
then squaring gives $(5\pm1)+(5\pm1)a^2+a(\pm1\pm1\pm1\pm1)=0$,
which is impossible. Similarly, no two of $v_i$ and $v_i\pm v_j$ ($i\ne j$)
are in the same line. So $\dim W_A=n-1>\frac12d(d-1)$.\qed

\begin{center}\textsc{4. Bounds on Universal Covers in $\R^3$}\end{center}

We give the proof of Theorem 3 as well as some non-optimal universal
covers, by truncating known universal covers in the literature
using the following simple lemma. The best upper bound for the minimum
volume of a universal cover in $\R^2$ was obtained by Gibbs and many
others before him using similar methods [Gi].

\lem (i) If a universal cover $K\subset\R^d$ has a reflectional symmetry
with normal $v\in S^{d-1}$, or a $\pi$-rotational symmetry around
a vector perpendicular to $v$, then $K\cap\{x\in\R^3:\<x,v\>\le\frac12\}$
is a universal cover.
\\(ii) If $K\subset K^d(v_1,\ldots,v_n)$ is a universal cover,
where $\|v_i\|=1$ for all $i$,
then the intersection of $K$ and the set of points at distance
$\le1$ to $K\cap F$ for every facet $F=\{y\in\R^d:\<y,v_i\>=\pm\frac12\}$
of $K^d(v_1,\ldots,v_n)$ is also a universal cover.

\pf (i) A set $X$ of diameter 1 cannot have points in both
$\{x\in\R^3:\<x,v\>>\frac12\}$ and $\{x\in\R^3:\<x,v\><-\frac12\}$,
so if $X$ is covered by $K$ then $X$ is also covered after removing one
of the sets. But due to symmetry, the result after removing each one
of the two sets from $K$ is congruent.

(ii) If a set $X$ of constant width 1 is covered by
$K\subset K^d(v_1,\ldots,v_n)$, then $X\cap F\ne\varnothing$
for every facet $F$, so $X$ does not contain points at distance $>1$ to 
any of the facets. Thus the set in the statement is a universal cover.\qed

\mn\textit{Decahedron (non-optimal).} By Proposition 5(i),
$K^3(v_1,\ldots,v_5)$ is a universal cover if $v_1,v_2,v_3$
are linearly independent and there is a rotational symmetry taking
$v_4$ to $-v_4$ and $v_5$ to $-v_5$. Gradient descent suggests
that the decahedron with minimum volume satisfying these conditions is
$K_1=K^3((0,0,1),(\cos\theta,0,\pm\sin\theta),(\cos\phi,\pm\sin\phi,0))$
for some $0<\theta\le\phi<\frac\pi2$. By Lemma 9(i), $K_1'=K_1\cap
\{x\in\R^3:x_1\le\frac12,x_2\cos\alpha+x_3\sin\alpha\le\frac12\}$
is a universal cover. The minimum volume of $K_1'$ is $\approx0.739534$,
where $(\theta,\phi,\alpha)\approx(0.640479,1.08474,0.78968)$.
By Lemma 9(ii), we obtain a smaller universal cover with volume
$\approx0.737006$.

\mn\textit{Truncation of octahedron (non-optimal).} Kuperberg proved in [K]
that $K_2=K^3(v_1^\pm,v_2^\pm,v_3^\pm)$ is a universal cover,
where $v_1^\pm=(0,\cos\theta,\pm\sin\theta)$,
$v_2^\pm=(\cos\theta,0,\pm\sin\theta)$, $v_3^\pm=(1,\pm1,0)$,
and $\theta\in(0,\frac\pi2)$. Using similar methods, Makeev proved in
[M6] that $K_3=K^3(v_1^\pm,v_2^\pm,e_1,e_2)$ is also a universal cover.
By Lemma 9(i),
$K_2'=K_2\cap\{x\in\R^3:\forall i,x_i\le\frac12\}$ and
$K_3'=K_3\cap\{x\in\R^3:x_3\le\frac12,\<x,v_3^\pm\>\le1/\sqrt2\}$
are universal covers. The volume of $K_2'$ is minimized when
$\theta=\frac\pi4$ (see below). The minimum volume of $K_3'$ is
$\approx0.691035$, where $\theta\approx0.823407$. By Lemma 9(ii),
we obtain a smaller universal cover with volume $\approx0.68483$.

\mn\textit{Proof of the upper bound in Theorem 3.} The dodecahedra
$K_2$ and $K_3$ above have symmetry groups $\Z/2\Z\times D_4$.
When $\theta=\frac\pi4$ in $K_2$, it becomes a rhombic dodecahedron
$R$ with octahedral symmetry $G=\Z/2\Z\times S_4$.

Let $\pm u_i\in S^2$ be vectors in the directions of the vertices
($1\le i\le7$), and $\pm w_i\in S^2$ vectors in the directions
of the perpendicular lines from 0 to the edges ($1\le i\le12$).
Let $X\subset R$ be a set of constant width 1, then for each $i$,
one of $X\cap\{x\in\R^3:\pm\<x,u_i\>>\frac12\}$ is empty, and similar
statements hold for $w_i$. This gives a map
$f:T=\{u_1,\ldots,u_7,w_1,\ldots,w_{12}\}\to\{\pm1\}$ of potential
truncations. Symmetries of $R$ permutes $T$. For a subset $T_0\subset T$,
if for any $f$ there exists $\gamma\in G$ such that
$\gamma T_0\subset f^{-1}(1)$, then
$R\cap\{x\in\R^3:\forall t\in T_0,\<x,t\>\le\frac12\}$
is a universal cover.

Brute-force search gives the five different maximal $T_0$ up to symmetry.
We shall give their description:
After scaling $R=K^3((0,1,\pm1),(\pm1,0,1),(1,\pm1,0))$, vectors
$u_1=e_1$, $u_2=e_2$, $u_3=e_3$ become vertices,
and let $u_4,u_5,u_6,u_7$ be the vertices connected to $\{e_1,e_2,e_3\}$,
$\{-e_1,e_2,e_3\}$, $\{e_1,-e_2,e_3\}$, $\{e_1,e_2,-e_3\}$ respectively.
We write $ij$ for the edge $(u_i,u_j)$ and $i\bar j$ for the edge
$(u_i,-u_j)$. Let $w_1,\ldots,w_{12}\in S^{12}$ be the
vectors in the directions of the perpendicular lines from 0 to the edges
34, 35, 36, $3\bar7$, 24, 25, $2\bar6$, 27, 14, $1\bar5$, 16, 17,
respectively. The five maximal $T_0$ which give universal covers
$R_1,\ldots,R_5$ are $\{u_1,u_2,u_3\}$, $\{u_1,u_2,u_5,u_6\}$,
$\{u_1,u_2,w_2\}$, $\{u_1,w_3,w_6\}$, $\{u_1,w_4,w_5\}$, respectively.
Their volumes are listed in the table below.

\begin{center}\begin{tabular}{lll}\hline
Set & Volume & $\approx$\\\hline
$R_1$&$\frac72-2\sqrt2$&0.671573\\$R_2$&
$\frac{11}{12}-\frac{23}{8}\sqrt2+\frac32\sqrt3+\frac12\sqrt6$&0.673624\\
$R_3$&$\frac73-\frac{15}8\sqrt2+\frac{13}{32}\sqrt6$&0.676788\\
$R_4,R_5$&$\frac76-\frac74\sqrt2+\frac{13}{16}\sqrt6$&0.682003\\
$R_1'$&&0.668948\\$R_2'$&&0.6698*\\$R_3'$&&0.6736*\\
$R_4',R_5'$&&0.6791*\\$R_1''$&&0.655984\\\hline
\end{tabular}\\ *Computed by triangle mesh which underestimates volume
\end{center}

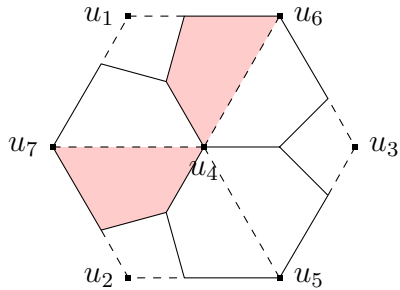
\begin{center}\begin{tikzpicture}[scale=2]
\coordinate (A) at (1,0); \coordinate (B) at (0.5,0.866);
\coordinate (C) at (-0.5,0.866); \coordinate (D) at (-1,0);
\coordinate (E) at (-0.5,-0.866); \coordinate (F) at (0.5,-0.866);
\coordinate (G) at (0.82,-0.32); \coordinate (H) at (0.5,0);
\coordinate (I) at (0.82,0.32); \coordinate (J) at (-0.13,0.866);
\coordinate (K) at (-0.25,0.433); \coordinate (L) at (-0.68,0.55);
\coordinate (M) at (-0.68,-0.55); \coordinate (N) at (-0.25,-0.433);
\coordinate (P) at (-0.13,-0.866);
\fill [color=red!20] (B) -- (J) -- (K) -- (0,0) -- cycle;
\fill [color=red!20] (D) -- (M) -- (N) -- (0,0) -- cycle;
\draw (K) -- (L) -- (D) -- (M) -- (N) -- (P) -- (F) -- (G);
\draw (G) -- (H) -- (I) -- (B) -- (J) -- (K) -- (0,0) -- (H);
\draw (N) -- (0,0); \draw [dashed] (J) -- (C) -- (L);
\draw [dashed] (M) -- (E) -- (P); \draw [dashed] (G) -- (A) -- (I);
\draw [dashed] (0,0) -- (B); \draw [dashed] (0,0) -- (D);
\draw [dashed] (0,0) -- (F); \Point{right}{$u_3$}{A};
\Point{left}{$u_1$}{C}; \Point{left}{$u_2$}{E}; \Point{left}{$u_7$}{D};
\Point{right}{$u_5$}{F}; \Point{right}{$u_6$}{B}; \Point{below}{$u_4$}{0,0}; 
\end{tikzpicture}\end{center}

The figure above shows three faces of $R_1$, which are parts of three
adjacent rhombic faces of $R$ after truncation at three vertices $u_1$,
$u_2$, $u_3$. Let $R_i'$ be the set obtained from $R_i$ after applying
Lemma 9(ii). $R_1'$ has $D_3$ symmetry around $u_4$. Since $X\subset R_1'$
is a set of constant width, it has at least one point on each of the three
faces that contain vertices $u_4,u_5,u_6,u_7$, which has been truncated at
$u_1,u_2,u_3$. By symmetry, without loss of generality, we may assume that
$X$ has points in the red regions in the figure, as well as points
in the corresponding red regions on the opposite faces.
Therefore, removing parts of $R_1'$ at distance $>1$ from each red region
as well as each corresponding red region on the opposite faces gives
a smaller universal cover $R_1''$ with volume $\approx0.655984$.\qed

\rem Y.\ Yang also found the universal cover $R_1'$ independently [Y].

\mn\textit{Proof of the lower bound in Theorem 3.} Actually we give
a lower bound for every dimension.

Any universal cover $K\subset\R^d$ covers a sphere $S$ of diameter 1
and a regular $d$-simplex $T$ of side length 1. Since $K$ is convex,
it contains $\mathrm{conv}(S\cup T)$, so a lower bound of the volume of
$K$ is the minimum volume of $\mathrm{conv}(S\cup T)$. 
Let $S$ be centered at 0, then the volume of
$\mathrm{conv}(S\cup\{(r,0,0)\})\setminus S$ is $f(r)=2^{-d}V_{d-1}
[\frac{2r}d(1-\frac1{4r^2})^{\frac{d+1}2}-
\int_{1/(2r)}^1(1-t^2)^{\frac{d-1}2}\,\d t]$ for $r\ge\frac12$,
and $f(r)=0$ for $r\le\frac12$, where $V_d$ is the volume of the unit ball
in $\R^d$. $f$ is increasing and $f''(r)>0$ for $r>\frac12$, so
$f$ is convex. Let $T=T_0+x$ where the regular $d$-simplex $T_0$
centered at 0 has vertices $v_0,\ldots,v_d$, then
$|\mathrm{conv}(S\cup T)|\ge2^{-d}V_d+\sum_{i=0}^df(\|v_i+x\|)=g(x)$.
Since $f$ is convex, $g$ is also convex, so it reaches minimum when
$\D g=0$ at $x=0$. In this case $|\mathrm{conv}(S\cup T)|=
2^{-d}V_d+(d+1)f((\frac d{2(d+1)})^{1/2})$. This is
$(\frac{5}{6\sqrt6}-\frac16)\pi\approx0.545193$ when $d=3$.\qed

\begin{center}\textsc{4. Universal Covers in Dimensions $d>2$}\end{center}

All cohomology in this section has $\F_2$-coefficients,
unless specified otherwise.

The proof of Theorem 4 uses the following lemma of Makeev.

\lem [M1] Let $K=K^d(v_1,\ldots,v_n)$ where $\|v_i\|=1$ for all $i$.
Let $G\subset O(d)$ be the group of symmetries of $K$,
and $V$ the $(n-d)$-dimensional $G$-representation on odd functions
on $\{\pm v_i\}_{i=1}^n$ modulo linear functions.
If there are no $G$-equivariant maps
$O(d)\to V\setminus\{0\}$, then $K$ is a universal cover.

\mn\textit{Proof of Theorem 4.}
We shall prove the case $d=2k$. The case $d=2k+1$ is similar.
The polygon $P_{2(4j-1)}$ has a central symmetry $c_j:(x,y)\mapsto(-x,-y)$
and a reflection symmetry $r_j:(x,y)\mapsto(x,-y)$, and
$\<c_j,r_j\>\cong\{\pm1\}^2$ acting as diagonal matrices.
We shall show that for our polytope, the conditions of Lemma 11 hold for
$G=\{\pm1\}^{2k}$ acting as diagonal matrices.

On the $(4j-1)$-dimensional space $W_j$ of odd functions on vertices of
$P_{2(4j-1)}$, $r_j$ has trace 1 (it fixes $(\pm1,0)$ and swaps other
opposite pairs of vertices with their reflections), so
$W_j\cong V_{j+}^{\oplus2j}\oplus V_{j-}^{\oplus(2j-1)}$,
where $c_j$ and $r_j$ act as $-1$ and
$\pm1$ respectively on $V_{j\pm}=\R$. On the 2-dimensional space $L_j$
of linear functions $\R^2\to\R$, $r_j$ has trace 0, so
$L_j\cong V_{j+}\oplus V_{j-}$. Therefore the $V$ in Lemma 11 is
$V=\bigoplus_{j=1}^k(W_j\ominus L_j)\cong
\bigoplus_{j=1}^k(V_{j+}^{\oplus(2j-1)}\oplus V_{j-}^{\oplus(2j-2)})$.

A nonvanishing $G$-equivariant map $O(d)\to V$ is equivalent to
a nonvanishing section of the real vector bundle
$O(d)\times_GV\to O(d)/G$. So it suffices to prove the Euler class
of this bundle is nonzero. Its classifying map factors through the
fiber bundle $O(d)/G\to BG\to BO(d)$. The map $BG\to BO(d)$ induces
$H^*(BO(d))\cong\F_2[\sigma_1,\ldots,\sigma_d]\to\F_2[x_1,\ldots,x_d]
\cong H^*(BG)$ where $\sigma_i$ is mapped to the $i$th elementary
symmetric polynomial. Using the Serre spectral sequence,
one may show that $H^*(O(d)/G)\cong
\F_2[x_1,\ldots,x_d]/(\sigma_1,\ldots,\sigma_d)$ [B1, Thm.~11.1],
where $x_j$ is the Euler class of the pullback of the tautological
bundle over the $j$th $B\{\pm1\}$. Thus $e(V_{j+})=x_{2j}$,
$e(V_{j-})=x_{2j-1}$, and
$e(V)=\prod_{j=1}^kx_{2j}^{2j-1}x_{2j-1}^{2j-2}=\prod_{j=1}^dx_j^{j-1}\ne0$
is the top class of $H^*(O(d)/G)$ [L, Ex.~4.19(b)].\qed

\rem Similar to Theorem 3, it is possible to improve the bounds on
the minimum volume of a universal cover in $\R^d$ by truncation
using Lemma 9. Similar to the argument in Theorem 3, the maximum number
of vertices one can truncate from the polytope $\prod_{j=1}^kP_{2n_j}$,
where $n_j$ is odd and $P_{2n_j}$ is the polygon with $2n_j$ vertices,
is a combinatorial problem: labeling the vertices by
$G=\prod_{j=1}^k\Z/2n_j\Z$ and the symmetries by
$T=\{(x_j)_{j=1}^k\mapsto(\eps_jx_j+a_j)_{j=1}^k:
\eps_j\in\{\pm1\},a_j\in\Z/2n_j\Z\}$, to find the largest
$S\subset G$ such that for every $f:G\to\{\pm1\}$
satisfying $f(x_1,\ldots,x_k)=-f(x_1+n_1,\ldots,x_k+n_k)$ for all
$(x_1,\ldots,x_k)\in G$, there exist $t\in T$ with $t(S)\subset f^{-1}(1)$.
A complicated argument, sketched in Appendix B, shows that the maximum 
$|S|$ is no more than 4 in all dimensions.

There are many other ways to truncate the polytope to get a smaller
universal cover. Similar to the argument in Theorem 3, one may truncate
a combination of vertices, edges, etc. One may also consider polytopes
with fewer facets but more symmetries than the one in Theorem 4, by
replacing the polygons $P_{2n_j}$ by ones with fewer edges, making
some of the $n_j$ equal, so that permuting them becomes additional symmetries.
One may also use products involving the rhombic dodecahedron $R$.
For example, $R^2$ is a universal cover in $\R^6$, and
$R\times P_{18}\times P_{26}$ is a universal cover in $\R^7$.
On the other hand, in high dimensions ($d>6$), even truncated polytopes
here are often bigger than Jung's ball, a ball of radius
$(\frac d{2(d+1)})^{1/2}$, which is an asymptotically optimal universal
cover [ABPR].

\medskip Finally, we prove Proposition 5 and Theorem 6.

\mn\textit{Proof of Proposition 5.} (i) When the condition is satisfied,
if there is a $\<g\>$-equivariant map $O(d)\to V\setminus\{0\}$, then
there is a path from 1 to $g$ in $SO(d)$ whose image in $V$ is a path
with endpoints $\pm1$. Taking quotient by $g$, we have a map
$SO(d)/\<g\>\to(V\setminus\{0\})/\{\pm1\}$ with nontrivial image
in $\pi_1$. But $SO(d)/\<g\>=\mathrm{Spin}(d)/\<-1,g\>$ has finite $\pi_1$
while $\pi_1((V\setminus\{0\})/\{\pm1\})=\Z$, a contradiction.
Now the result follows from Lemma 11.

(ii) This is known when $d$ is odd (see Proposition 14),
so we assume $d$ is even. By Lemma 11, it suffices to prove that
there are no continuous odd maps $f:SO(d)\to S^k$.
We also denote by $f$ the induced map $PSO(d)\to\RP^k$.
We have Gysin sequences
\begin{center}\begin{tikzcd}[row sep=scriptsize,column sep=small]
H^{k-1}(\RP^k)\ar[r,"w"]\ar[d,"f^*"]&H^k(\RP^k)\ar[r]\ar[d,"f^*"]&
H^k(S^k)\ar[r,"i"]\ar[d,"f^*"]&H^k(\RP^k)\ar[r]\ar[d,"f^*"]&0\ar[d]\\
H^{k-1}(PSO(d))\ar[r,"f^*w"]&H^k(PSO(d))\ar[r]&H^k(SO(d))\ar[r]&
H^k(PSO(d))\ar[r]&H^{k+1}(PSO(d))\end{tikzcd}\end{center}
where in $H^*(\RP^k)\cong\F_2[w]/(w^{k+1})$ we have an isomorphism
$H^{k-1}(\RP^k)\to H^k(\RP^k)$,
so $i:H^k(S^k)\to H^k(\RP^k)$ is an isomorphism, taking the generator $u$
to $w^k$. The first Stiefel--Whitney class $f^*w$ of $SO(d)\to PSO(d)$
satisfies $f^*w^{2^e-1}\ne0$ [B2, Thm.~11.5].

Note that in the case $d=2^e$ we already get a contradiction:
$f^*w^{k+1}\ne0$ but $w^{k+1}=0$. Now we assume that $d\ne2^e$.
In \[ H^*(SO(d))\cong\F_2[b_i:1\le i<d,i\text{ odd}]/(b_i^{p_i}), \]
where $\deg b_i=i$ and $p_i$ is the smallest power of 2 no less than $d/i$, 
$f^*u\ne0$ is a sum of monomials of the form $\prod_ib_i^{e_i}$.
Since $f^*u^2=0$ is the sum of monomials $\prod_ib_i^{2e_i}$,
we have $2e_i\ge p_i$ for some $i$,
so $d>2k\ge2ie_i\ge ip_i\ge d$, a contradiction.

When $e=1$, the argument above shows there are no continuous odd maps
$SO(d)\to S^1$. Conner and Floyd showed that, if there is a continuous odd
map $f:SO(d)\to S^2$, then $f^*:H^2(S^2;\Z)\to H^2(SO(d);\Z)$ is
nontrivial [CF, Thm.~4.6] (the condition of their theorem is that the
primary obstruction to an odd continuous map $SO(d)\to S^1$ is nonzero;
but that is the only obstruction here).

It is not hard to show by induction that $H_1(SO(d);\Z)\cong\Z/2\Z$
and $H_2(SO(d);\Z)=0$ for $d\ge3$ using the Serre spectral sequence for
the fiber bundle $SO(d-1)\to SO(d)\to S^{d-1}$,
starting from $SO(3)\cong\RP^3$. By the universal coefficient theorem,
$H^1(SO(d);\Z)=0$ and $H^2(SO(d);\Z)\cong\Z/2\Z$, so the integral Bockstein
$\beta:H^1(SO(d);\F_2)\to H^2(SO(d);\Z)$ is an isomorphism.
Let $u$ generate $H^2(S^2;\Z)$, then by the previous paragraph,
$f^*u\ne0$ in $H^2(SO(d);\Z)$, so $f^*u=\beta(b_1)$.
Let $\rho:H^2(SO(d);\Z)\to H^2(SO(d);\F_2)$ be the reduction mod 2,
then $\rho\beta(x)=x^2$ (it is $\mathrm{Sq}^1$). By the structure
of $H^*(SO(d);\F_2)$, for $d\ge6$, we have
$0\ne b_1^4=\rho\beta(b_1^2)=\rho(f^*u^2)=0$, a contradiction.
So there are no continuous odd maps $SO(d)\to S^2$ when $e=1$.

Makeev gave the following simpler proof for $d=2^e$, $e\in\{2,3\}$ in [M1]:
We have odd maps $S^{d-1}\to SO(d)$ given by quaternion and
octonion multiplication, so there are no odd maps $SO(d)\to S^{d-2}$
by the Borsuk--Ulam theorem.\qed

\mn\textit{Proof of Theorem 6.} We only prove the case when $d=2k+1$;
the even case is similar. We have $K=D\times P^{k-1}=K^d(v_1,\ldots,v_{2d})$,
where $\{v_1,\ldots,v_{2d}\}=\{e_1,\ldots,e_d\}
\cup\{(e_1+e_2)/\sqrt2,(e_2+e_3)/\sqrt2,(e_3+e_1)/\sqrt2\}\cup
\{(e_{2j}\pm e_{2j+1})/\sqrt2:2\le j\le k\}$.
By Lemma 8, there is a set $X\subset\R^d$
of constant width 1 such that $s_X(v)=\eps\<e_1,v\>^3$ for $v\in S^{d-1}$.
Suppose $K$ is a universal cover, then in particular
there exist $A\in O(d)$ and $b\in\R^d$ such that $s_X(Av_i)=\<b,Av_i\>$ for
$1\le i\le2d$.

Let $c=\eps^{-1}A^\top b$ and $a=A^\top e_1$, then
$\<c,v_i\>=\eps^{-1}\<b,Av_i\>=\<e_1,Av_i\>^3=\<a,v_i\>^3$.
In particular, $c_i=a_i^3$, so
$2(a_{2j}^3\pm a_{2j+1}^3)=(a_{2j}\pm a_{2j+1})^3$ for $2\le j\le k$.
By adding and subtracting these two equations, we have
$a_{2j}(a_{2j}^2-3a_{2j+1}^2)=a_{2j+1}(a_{2j+1}^2-3a_{2j}^2)=0$,
which implies $a_{2j}=a_{2j+1}=0$. Similarly, we have the
three equations $2(a_i^3+a_j^3)=(a_i+a_j)^3$ for $(i,j)=(1,2),(2,3),(3,1)$.
If one of $a_1,a_2,a_3$ is 0, then the equations imply that
the other two are 0, which is impossible because $\|a\|=1$.
Otherwise, by the three equations,
$a_1/a_2,a_2/a_3,a_3/a_1\in\{-1,2\pm\sqrt3\}$, which is also impossible
because their product cannot be 1.\qed

\begin{center}\textsc{Appendix A. Some Known Results}\end{center}

We reproduce Makeev's proof of the upper bound on the number
of facets of the polytope universal cover for the reader's convenience.

\prop [M2] No $K^d(v_1,\ldots,v_n)$ with $n>\frac12d(d+1)$
is a universal cover.

\pf Let $\Omega$ be the set of odd smooth functions on $S^{d-1}$,
and $f:O(d)\times\Omega\to\R^{n-d}$ the map $(A,s)\mapsto
(s(Av_i)-\<b,Av_i\>)_{i=d+1}^n$, where $s(Av_i)=\<b,Av_i\>$ for
$1\le i\le d$.

For $d<i\le n$ and $A\in O(d)$, let $s_i^A\in\Omega$ such that
$s_i^A=\delta_{ij}$ near $Av_j$, then
$\forall i,s_i^A(A'v_j)=\delta_{ij}$ for $A'$ in an open neighborhood
$U_A\subset O(d)$ of $A$. Let $O(d)\subset\bigcup_kU_{A_k}$ be a finite
subcover, $s_0\in\Omega$ a given function, and
$V=\mathrm{span}\{s_0,s_i^{A_k}\}\subset\Omega$.

For a given point $(A,s)\in V$, let $A\subset U_{A_k}$, then
the restriction of $f_0=f|_{O(d)\times V}$ to
$(A,s+\sum_{i=d+1}^nt_is_i^{A_k})$ for $t\in\R^{n-d}$ has derivative
the identity at $t=0$. So $f_0$ is a submersion and 
$\dim f_0^{-1}(0)<\dim V$. Thus the projection of
$f_0^{-1}(0)$ to $V$ has measure 0. Therefore the set of
$s_0\in\Omega$ such that $f(A,s)\ne0$ for all $A\in O(d)$
is a dense open set.\qed

\medskip When $d$ is odd, Makeev proved the following special case of
Proposition 5(ii).

\prop [M1] Any $K^d(v_1,\ldots,v_{d+1})$ is a universal cover.

\pf Given a set $X$ of constant width 1, we need to find $A,b$
such that $s_X(Av_i)=\<b,Av_i\>$ for $1\le i\le n$. Without loss of
generality, $v_1,\ldots,v_d$ are linearly independent.
Let $v_{d+1}=\sum_{i=1}^da_iv_i$, then $b$ exists iff
$s_X(Av_{d+1})-\sum_{i=1}^da_is_X(Av_i)=0$. For each $i$, there
is a rotation taking $v_i$ to $-v_i$, so $\int_{SO(d)}s_X(Av_i)\,\d A=0$.
Thus $\int_{SO(d)}[s_X(Av_{d+1})-\sum_{i=1}^da_is_X(Av_i)]\,\d A=0$,
and the integrand cannot be nonzero.\qed

\begin{center}\textsc{Appendix B. Truncation of the Polytope in
Theorem 4}\end{center}

Remark 12 reduces the problem of maximum number of vertices that
can be truncated using Lemma 9 from the polytope in Theorem 4 to
the combinatorial Theorem 16 below.

\prop Let $n>1$ be an odd integer and $S\subset\Z/n\Z$. The following
are equivalent:\\(i) for every $f:\Z/2n\Z\to\{\pm1\}$, there exist
$(\eps,a)\in\{\pm1\}\times\Z/n\Z$ such that $|f(\eps S+a)|=1$.
(ii) Either $|S|\le2$ or there exist
$(\eps,a)\in\{\pm1\}\times\Z/n\Z$ such that
$\eps S+a=\{0,\frac n7,\frac{3n}7\}$.

\mn\textit{Proof Sketch.} The only hard case is to prove (i) implies (ii)
when $|S|=3$. Then the $|S|>3$ case also holds as subsets of $S$
also satisfy (i). Without loss of generality, $S=\{0,x,y\}$ where
$0<x<y<\frac n2$, and $\gcd(x,y,z,n)=1$ (by dividing).
Let $z=y-x$, and let $\<a\>\in(-\frac n2,\frac n2)$
such that $\<a\>\equiv a\mod n$. 

By considering the coloring $f(x)=1\iff\<kx\>>0$, we know that the pairwise
differences $\<kx\>,-\<ky\>,\<kz\>\in(-\frac n2,\frac n2)$ do not have
the same sign. They sum to a multiple of $n$, so the sum must be 0.
Let $f_x(k)=\<kx\>$, then $f_x+f_z=f_y$.

Consider the discrete Fourier transform $\hat f(j)=\sum_{k\in\Z/n\Z}
f(k)e^{-2\pi ijk/n}$. We have $\hat f_x(j)=0$ when $d=\gcd(x,n)\nmid j$,
and $\hat f_x(j)=d^2g_{n_0}(x_0^{-1}j_0)$ otherwise, where
$x_0=x/d$, $n_0=n/d$, $j_0=j/d$, and $g_m(j)=\frac{(-1)^jim}{2\sin(\pi j/m)}$
for $j\in\Z/m\Z$. We have $\hat f_x+\hat f_z=\hat f_y$.
A simple casework using this formula shows $\gcd(x,n)=\gcd(y,n)=\gcd(z,n)=1$.

We have $\<(k+1)x\>-\<kx\>=x-k1_{\<kx\>>n/2-x}$, so $\<kx\>+\<kz\>=\<ky\>$
and $x+z=y$ imply $K_x\sqcup K_z=K_y$ where
$K_x:=\{k\in\Z/n\Z:\<kx\>>\frac n2-x\}=x^{-1}(\frac n2-x,\frac n2)$.
We may compute $12\sum_{k\in K_x}k^2=4x-x^{-1}$ and
$240\sum_{k\in K_x}k^4=48x-40x^{-1}+7x^{-3}$, so $x^{-1}+z^{-1}=y^{-1}$
and $7(x^{-3}+z^{-3}-y^{-3})=0$. Cubing, we have
$y^{-3}=x^{-3}+z^{-3}+3(x^{-1}+z^{-1})(xz)^{-1}
=x^{-3}+z^{-3}+3(xyz)^{-1}$, so $21(xyz)^{-1}=0$. Hence $n\mid 21$.
Since $0<x<y<\frac n2$, we have $n\in\{7,21\}$.

If $n=21$, since $x^{-1}+z^{-1}=y^{-1}$, we have $x^2-xy+y^2=0$.
Thus $3\mid x+y$, and $7S=\{0,7,14\}$, so $S$ fails the condition (i).
If $n=7$, it is not hard to check the result.\qed

\thm Let $n_j$ be odd integers, $G:=\prod_{j=1}^k\Z/2n_j\Z$,
$(x_j)^*:=(x_j+n_j)$, $T:=\{G\to G,(x_j)\mapsto(\eps_jx_j+a_j):
\eps_j\in\{\pm1\},a_j\in\Z/2n_j\Z\}$, and $S\subset G$. If for every
$f:G\to\{\pm1\}$ satisfying $f(x^*)=-f(x)$ for all $x\in G$, there
exists $t\in T$ such that $t(S)\subset f^{-1}(1)$, then $|S|\le4$.

\mn\textit{Proof Sketch.} Suppose $S$ satisfy the condition.
Then it is in one coset of $\prod_{j=1}^k\Z/n_j\Z\subset G$
(otherwise it has both even and odd projections in coordinate $j$,
and thus fails the condition for $f(x)=(-1)^{x_j}$).
So it suffices to consider one coset: we replace
$G$ with $\prod_{j=1}^k\Z/n_j\Z$ where $n_j>1$ are odd, and the condition
is that for every $f:G\to\{\pm1\}$ there exists $t\in T$ such that
$|f(t(S))|=1$.

Suppose $|S|=5$. The projection $p_j(S)$ of $S$ to each coordinate $j$
satisfies the condition of Proposition 15. By considering further cosets
of subgroups of $G$, we may assume either $p_j(S)=\{0,1\}$ or
$p_j(S)=\{0,1,3\}$ with $n_j=7$. Let $J$ be a minimal
set of coordinates to which the projection of $S$ is injective,
then $|J|<|S|=5$. Up to translations in $T$, there are two orderings of
$\{0,1,3\}$ and one ordering of $\{0,1\}$.

Furthermore, for those coordinates $j$ where $p_j(S)=\{0,1\}$ and
$n_j>5$, it suffices to find a coloring $f$ with $n_j=5$ where $S$
fails the condition. Then we get a coloring for the original $G$
by composing $f$ and the map taking $0,\ldots,n_j-1$ to
$0,1,2,3,4,3,4,\ldots,3,4$.

These give finitely many possible projections of $S$ to $|J|\le4$ coordinates.
A computer program found a coloring for each of them making $S$ fail
the condition.\qed

\begin{center}\textsc{References}\end{center}\parindent=0pt

\let\pbr=\par\def\par{\hangindent=2em\pbr}

[ABPR] A.\ Arman, A.\ Bondarenko, A.\ Prymak, and D.\ Radchenko (2025),
On asymptotic Lebesgue’s universal covering problem, arXiv:2512.04023.

[B1] A.\ Borel (1953), La cohomologie mod 2 de certains espaces homog\`enes,
\textsl{Comment.\ Math.\ Helv.} \textbf{27}, 165--197.

[B2] ------ (1954), Sur l'homologie et la cohomologie des groupes de
Lie compacts connexes, \textsl{Am.\ J.\ Math.} \textbf{76}:2, 273--342.

[BS] P.\ Brass and M.\ Sharifi (2005), A lower bound for Lebesgue's
universal cover problem, \textsl{Int.\ J.\ Comput.\ Geom.\ Appl.}
\textbf{15}, 537--544.

[C] H.\ S.\ M.\ Coxeter (1988), Regular and semi-regular polytopes. III.
\textsl{Math.\ Z.} \textbf{200}, 3--45.

[CF] P.\ E.\ Conner and E.\ E.\ Floyd (1960), Fixed point free involutions
and equivariant maps, \textsl{Bull.\ Am.\ Math.\ Soc.}
\textbf{66}(6): 416--441.

[CS] J.\ H.\ Conway and N.\ J.\ A.\ Sloane (1999), \textsl{Sphere Packings,
Lattices and Groups}, 3rd ed.

[DGS] P.\ Delsarte, J.\ M.\ Goethals and J.\ J.\ Seidel (1977),
Spherical codes and designs, \textsl{Geom.\ Dedicata} \textbf6, 363--388.

[Gi] P.\ Gibbs (2018), An upper bound for Lebesgue's covering problem,
arXiv:1810.10089.

[Gr] P.\ M.\ Gruber (2007), \textsl{Convex and Discrete Geometry}.

[GS] J.\ M.\ Goethals and J.\ J.\ Seidel (1979), Spherical designs.
\textsl{Proc.\ Sympos.\ Pure Math.} \textbf{34}, 255--272.

[K] G.\ Kuperberg (1999), Circumscribing constant-width bodies with
polytopes, \textsl{New York J.\ Math.} \textbf5, 91--100.

[M1] V.\ V.\ Makeev (1981), Universal'nyye pokryshki. I [Universal covers. I], \textsl{Ukr.\ Geom.\ Sb.} \textbf{24}, 70--79 (in Russian). (Original Cyrillic: В.\ В.\ Макеев, Универсальные покрышки, \textsl{Укр.\ Геом.\ Сб.})

[M2] ------ (1982), Universal'nyye pokryshki. II [Universal covers. II],
\textsl{Ukr.\ Geom.\ Sb.} \textbf{25}, 82--86 (in Russian). (Original Cyrillic: В.\ В.\ Макеев, Универсальные покрышки. \textsl{Укр.\ Геом.\ Сб.})

[M3] ------ (1982), Razmernostnyye ogranicheniya v zadachakh kombinatornoy geometrii [Dimensional restrictions in problems of combinatorial geometry], \textsl{Sibirsk. Mat. Zh.} \textbf{23}:4, 197--201 (in Russian). (Original Cyrillic: В.\ В.\ Макеев, Размерностные ограничения в задачах комбинаторной геометрии, \textsl{Сиб.\ матем.\ журн.})

[M4] ------ (1984), Primeneniye topologii v nekotorykh zadachakh kombinatornoy geometrii [Applications of topology to certain problems in combinatorial geometry], \textsl{Ukr.\ Geom.\ Sb.} \textbf{27}, 83--88 (in Russian). (Original Cyrillic: В.\ В.\ Макеев, Применение топологии в некоторых задачах
комбинаторной геометрии, \textsl{Укр.\ Геом.\ Сб.}) 

[M5] ------ (1997), Of affine images of a rhombododecaedron circumscribed about a convex body in $\R^3$, in \textsl{Geometry and Topology. Part 2}, Zap.\ Nauchn.\ Sem.\ POMI, \textbf{246}, POMI, St.\ Petersburg, 191--195; \textsl{J.\ Math.\ Sci. (New York)}, \textbf{100}:3 (2000), 2307--2309.

[M6] ------ (2000), Three-dimensional polytopes inscribed in and circumscribed about convex compacta, \textsl{Algebra i Analiz}, \textbf{12}:4, 1--15; \textsl{St.\ Petersburg Math.\ J.}, \textbf{12}:4 (2001), 507--518.

[MMO] H.\ Martini, L.\ Montejano, and D.\ Oliveros (2019),
\textsl{Bodies of Constant Width}.

[NS] M.\ V.\ Noskov and H.\ J.\ Schmid (2004), On the number of nodes
in $n$-dimensional cubature formulae of degree 5 for integrals
over the ball, \textsl{J.\ Comput.\ Appl.\ Math.} \textbf{169}:2, 247--254.

[Y] Y.\ Yang (2026), Further Sprague-Type Truncations of the Makeev
Polyhedron for the Three-Dimensional Lebesgue Universal Covering Problem.
Preprint,\\doi: 10.5281/zenodo.21349947.
\end{document}